%% file: h2bf2r_arxiv.tex
\documentclass[a4paper,12pt]{article}

\usepackage{amsmath,amssymb,amsthm}

\newcommand{\R}{\mathbb{R}}
\newcommand{\Z}{\mathbb{Z}}
\newcommand{\N}{\mathbb{N}}
\newcommand{\Hom}{\mathrm{Hom}}
\newcommand{\PSL}{\mathrm{PSL}}
\newcommand{\deff}{\mathrm{def}\,}
\newcommand{\h}{\mathrm{H}}
\newcommand{\eh}{\mathrm{EH}}
\newcommand{\bb}{\mathrm{b}}
\newcommand{\hh}{\mathrm{h}}
\newcommand{\im}{\mathrm{im}\,}
\newcommand{\QM}{\mathrm{QM}}

\begin{document}

\title{Quasi-morphisms on Free Groups}
\author{Pascal Rolli}
\date{}
\maketitle

\vspace{1cm}

\begin{abstract}
Let $F$ be the free group over a set of two or more generators. In \cite{Broo} R. Brooks constructed an infinite family of quasi-morphisms $F\rightarrow\R$ such that an infinite subfamily gives rise to independent classes in the second bounded cohomology $\h^2_\bb(F,\R)$, which proves that this space is infinite dimensional, cf. \cite{Mit}. We give a simpler proof of this fact using a different type of quasi-morphisms. After computing the Gromov norm of the corresponding bounded classes, we generalize our example to obtain quasi-morphisms on free products, as well as quasi-morphisms into groups without small subgroups, also known as $\varepsilon$-representations.\\\\
This article is an extract from the author's Master thesis, conducted under the supervision of Alessandra Iozzi at ETH Z\"urich.
\end{abstract}

\vspace{1cm}

\tableofcontents

\newpage

\input{section1.tex}
\input{section2.tex}

\input{section3.tex}
\input{section4.tex}
\input{section5.tex}

\newpage

\bibliographystyle{acm}
\bibliography{bibliography}
\vspace{2cm}
Department of Mathematics\\
ETH Z\"urich, CH-8092 Z\"urich, Switzerland\\\\
\texttt{prolli@student.ethz.ch}

\end{document}

%% file: section1.tex
\section{Introduction}
Let $\Gamma$ be a group. A map $f:\Gamma\rightarrow\R$ is called a \textit{quasi-morphism} if
\[
\sup_{x,y\in\Gamma}|f(x)+f(y)-f(xy)|<\infty
\]
If $\varphi:\Gamma\rightarrow\R$ is a homomorphism and $b:\Gamma\rightarrow\R$ is a bounded map then $f=\varphi+b$ is a quasi-morphism. Any quasi-morphism of this form is said to be \textit{trivial}. Given a group, how many (in a suitable sense) non-trivial quasi-morphisms does it admit? An answer to this question may be found in the second bounded cohomology.\\\\
Recall that the ordinary group cohomology of $\Gamma$ with (trivial) real coefficients, denoted $\h^*(\Gamma,\R)$, can be defined as the homology of the complex $\left(C^*(\Gamma),\,\partial\right)$ where $C^n(\Gamma)$ is the space of maps $\Gamma^n\rightarrow\R$ and $\partial:C^n(\Gamma)\rightarrow C^{n+1}(\Gamma)$ is given by
\begin{align*}
\partial f(x_0,\dots,x_{n})&=f(x_1,\dots,x_{n})+\sum_{i=0}^{n-1}(-1)^if(x_0,\dots,x_ix_{i+1},\dots,x_{n})\\&+(-1)^nf(x_0,\dots,x_{n-1}).
\end{align*}
Restricting the boundary maps $\partial$ to the subspaces
\[
C^n_\bb(\Gamma):=\{f:\Gamma\rightarrow\R\,;\,f\mbox{ is bounded}\}\subset C^n(\Gamma)
\]
gives a subcomplex $(C^*_\bb(\Gamma),\partial)$ whose homology $\h^*_\bb(\Gamma,\R)$ is called the \textit{bounded cohomology} of $\Gamma$. More precisely, let
\begin{align*}
Z^n_\bb(\Gamma):=\ker\left(\partial:C^n_\bb(\Gamma)\rightarrow C^{n+1}_\bb(\Gamma)\right)\\
B^n_\bb(\Gamma):=\im\left(\partial:C^{n-1}_\bb(\Gamma)\rightarrow C^{n}_\bb(\Gamma)\right)
\end{align*}
be the spaces of \textit{bounded $n$-cocycles} and \textit{bounded $n$-coboundaries} respectively. The $n$-th bounded cohomology of $\Gamma$ is then the quotient
\[
\h^n_\bb(\Gamma,\R):=\frac{Z^n_\bb(\Gamma)}{B^n_\bb(\Gamma)}.
\]
The inclusion $C^*_\bb(\Gamma)\hookrightarrow C^*(\Gamma)$ induces the \textit{comparison map}
\[
c:\h^*_\bb(\Gamma,\R)\rightarrow\h^*(\Gamma,\R)
\]
Its kernel in degree two, denoted $\eh^2_\bb(\Gamma,\R)$, turns out to be relevant to the study of quasi-morphisms:\\\\
Let $\QM(\Gamma)\subset C^1(\Gamma)$ denote the space of quasi-morphisms on $\Gamma$. For $f\in\QM(\Gamma)$, the quantity $\sup_{x,y\in\Gamma}|f(x)+f(y)-f(xy)|=:\deff f$ is called the \textit{defect} of $f$. Since $\partial f(x,y)=f(y)-f(xy)+f(x)$ we get
\[
\|\partial f\|_\infty=\deff f<\infty,
\]
so $\partial f$ is a bounded 2-cocycle and we have a linear map
\[
\Phi: \QM(\Gamma)\rightarrow \h^2_\bb(\Gamma,\R),\quad f\mapsto[\partial f]_\bb,
\]
where $[\,\cdot\,]_\bb$ denotes the bounded cohomology class.\\\\
One has $[\partial f]_\bb=0$ if and only if there is some $b\in C^1_\bb(\Gamma)$ such that $\partial(f-b)=0$, or equivalently, such that $f-b$ is a homomorphism. So the kernel of the above map is the subspace of trivial quasi-morphisms:
\[
\ker\Phi=\Hom(\Gamma,\R)\oplus C^1_\bb(\Gamma)\subset\QM(\Gamma).
\]
Note that the class $[\partial f]_\bb$ lies in the kernel of the comparison map, since the corresponding unbounded class $[\partial f]\in\h^2(\Gamma,\R)$ is zero. On the other hand, if $\alpha\in\h^2_\bb(\Gamma,\R)$ satisfies $c(\alpha)=0$ then $\alpha=[\partial f]_\bb$ for some $f\in C^1(\Gamma)$ which must be a quasi-morphism. That is, $\im\Phi=\eh^2_\bb(\Gamma,\R)$, so we have the following\\\\
\textbf{Proposition 1.1.} \textit{The map $\Phi$ induces an isomorphism}
\[
\frac{\QM(\Gamma)}{\Hom(\Gamma,\R)\oplus C^1_\bb(\Gamma)}\cong\eh^2_\bb(\Gamma,\R).
\]
\smallskip
\textbf{Corollary 1.2.} \textit{If $\h^2(\Gamma,\R)=0$ then $\Phi$ induces an isomorphism}
\[
\frac{\QM(\Gamma)}{\Hom(\Gamma,\R)\oplus C^1_\bb(\Gamma)}\cong\h^2_\bb(\Gamma,\R).
\]

%% file: section2.tex
\section{The second bounded cohomology of a free group}
Let $F=F(S)$ be the free group over a set $S$, $|S|\geq2$. Since $\h^2(F,\R)=0$ (cf. \cite{Brow}), the calculation of $\h^2_\bb(F,\R)$ amounts to finding non-trivial quasi-morphisms on $F$. In \cite{Broo}, R. Brooks constructed an infinite family of quasi-morphisms on $F$ such that an infinite subfamily is mapped to independent classes under $\Phi$ (cf. \cite{Mit}), which proves that space $\h^2_\bb(F,\R)$ is infinite dimensional. In the following we describe another type of quasi-morphisms on $F$ which allows us to give a simpler proof of this fact.\\\\
We say $x\in F$ is a \textit{power} if $x=s^k$ for some $s\in S$ and some $k\in\Z$. Each non-trivial
$x\in F$ has a unique shortest factorization into powers, which we simply call \textit{factorization of} $x$. Let $\ell^\infty$ be the space of bounded real sequences. For $\sigma\in\ell^\infty$ we define a map $g_\sigma:F\rightarrow\R$ as follows: For a power $x=s^k$ set $g_\sigma(x)=\sigma(k)$, where $\sigma$ is extended to an odd function on $\Z$. In general, for $x\in F$ with factorization $x=x_0\cdots x_n$, set
\[
g_\sigma(x)=\sum_{i=0}^n g_\sigma(x_i).
\]
\textbf{Proposition 2.1.} \textit{The map $g_\sigma$ is a quasi-morphism.}\\\\
\textbf{Proof.} Let $x,y\in F$ have factorizations $x=x_0\cdots x_n$ and $y=y_0\cdots y_m$. The factorization of $xy$ has the form
\[x_0\cdots x_{n-r}\cdot z\cdot y_r\cdots y_m\qquad\mbox{or}\qquad x_0\cdots x_{n-r}\cdot y_r\cdots y_m\]
for some $r\geq0$. The first case occurs if $x_{n-i}=y_i^{-1}$ for $0\leq i\leq r-2$. Since $\sigma$ is odd, this implies $g_\sigma(x_{n-i})+g_\sigma(y_i)=0$, and hence,
\begin{align*}
\left|g_\sigma(x)+g_\sigma(y)-g_\sigma(xy)\right|
&=\left|\sum_{i=0}^{n-r+1}g_\sigma(x_i)+\sum_{i=r-1}^m g_\sigma(y_i)-g_\sigma(xy)\right|\\
&=\left|g_\sigma(x_{n-r+1})+g_\sigma(y_{r-1})-g_\sigma(z)\right|\\
&\leq3\left\|\sigma\right\|_\infty.
\end{align*}
In the second case $x_{n-i}=y_i^{-1}$ holds as well for $i=r-1$, so $\partial g_\sigma(x,y)=0$.\qed\\\\
\textbf{Proposition 2.2.} \textit{The linear map $\ell^\infty\rightarrow\h^2_\bb(F,\R)$, $\sigma\mapsto[\partial g_\sigma]_\bb$ is injective.}\\\\
\textbf{Proof.} Assume that $[\partial g_\sigma]_\bb=0$. This means $g_\sigma\in\ker\Phi$, i.e. $g_\sigma-b=\varphi$ for some $b\in C^1_\bb(F)$ and some $\varphi\in\Hom(F,\R)$. For $s\in S$, evaluating this equation at $s^k$ yields $\sigma(k)-b(s^k)=k\,\varphi(s)$. The left-hand side is bounded as a function of $k$, so $\varphi(s)=0$. Hence, $\varphi=0$ and $g_\sigma$ is bounded.\\
Let $s,t\in S$ be two distinct generators. For $k,l\in\Z$ the equation $g_\sigma((s^lt^l)^k)=2k\,\sigma(l)$ holds. Since $g_\sigma$ is bounded, this implies $\sigma(l)=0$, and so $\sigma=0$.\qed\\\\
\textbf{Corollary 2.3.} \textit{The space $\h^2_\bb(F,\R)$ has infinite dimension.}\\\\
\textbf{Remarks.} (i) Note that the argument holds with slight modification if we define $g_\sigma(s^k)=\sigma_s(k)$ where $\sigma=(\sigma_s)_{s\in S}\in\left(\ell^\infty\right)^S$ is a uniformly bounded family of sequences, cf. Section 4.\\\\
(ii) Free groups $F(S)$ as above belong to the class of \textit{non-elementary hyperbolic} groups. D.B.A. Epstein and K. Fujiwara proved that $\h^2_\bb(\Gamma,\R)$ is infinite dimensional for any such group $\Gamma$. (\cite{EpsFuj}).\\\\
(iii) For a group $\Gamma$, a complex Hilbert space $\mathcal{H}$ and a unitary representation $\pi:\Gamma\rightarrow U(\mathcal{H})$, the bounded cohomology of $\Gamma$ \textit{with coefficients in $\mathcal{H}$}, denoted $\h^*_\bb(\Gamma,\mathcal{H})$, is defined as the homology of the complex $(C^*_\bb(\Gamma,\mathcal{H}),\partial_\pi)$. Here $C^n_\bb(\Gamma,\mathcal{H})$ is the space of bounded maps $\Gamma^n\rightarrow\mathcal{H}$ and $\partial_\pi:C^n_\bb(\Gamma,\mathcal{H})\rightarrow C^{n+1}_\bb(\Gamma,\mathcal{H})$ is given by
\begin{align*}
\partial_\pi f(x_0,\dots,x_{n})&=\pi(x_0)f(x_1,\dots,x_{n})+\sum_{i=0}^{n-1}(-1)^if(x_0,\dots,x_ix_{i+1},\dots,x_{n})\\&+(-1)^nf(x_0,\dots,x_{n-1}).
\end{align*}
It seems to be unknown whether $\h^2_\bb(F,\mathcal{H})\neq0$ for non-trivial representations $\pi:F\rightarrow U(\mathcal{H})$.\\
Consider the following generalization of the above construction: Let $\sigma=(\sigma_s)_{s\in S}\in\left(\ell^\infty(\mathcal{H})\right)^S$ be a uniformly bounded family of sequences in $\mathcal{H}$, where each $\sigma_s$ is extended to a map on $\Z$ such that $\sigma_s(k)+\pi(s^k)\sigma_s(-k)=0$ for $k\in\Z$. Define $g_\sigma:F\rightarrow\mathcal{H}$ on powers by $g_\sigma(s^k)=\sigma_s(k)$, and for $x\in F$ with factorization $x=x_0\cdots x_n$ set
\begin{align*}
g_\sigma(x)&:=\\&g_\sigma(x_0)+\pi(x_0)g_\sigma(x_1)+\pi(x_0x_1)g_\sigma(x_2)+\dots+\pi(x_0x_1\cdots x_{n-1})g_\sigma(x_n).
\end{align*}
As above $\|\partial_\pi g_\sigma\|_\infty<\infty$, so there is a linear map $\left(\ell^\infty(\mathcal{H})\right)^S\rightarrow\eh^2_\bb(F,\mathcal{H})$, $\sigma\mapsto[\partial_\pi g_\sigma]_\bb$. The space $\h^2_\bb(F,\mathcal{H})$ could be shown to be non-zero by proving that the image of some $\sigma\in\left(\ell^\infty(\mathcal{H})\right)^S$ is a non-trivial bounded class, which means that $g_\sigma$ is not the sum of a crossed homomorphism and a bounded map. This is most likely to work if $\pi$ has a trivial subrepresentation.
\\\\
We conclude this section with the following observation (which concerns the case of trivial real coefficients):\\\\
\textbf{Proposition 2.4.} $\|\partial g_\sigma\|_\infty=\deff g_\sigma=\deff\sigma$.\\\\
\textbf{Proof.} If $x,y\in F$ are such that the first case occurs in the proof of Proposition 2.1, then in fact the three terms of $\partial g_\sigma(x,y)$ that remain after cancellation are powers of the same generator. So $|\partial g_\sigma(x,y)|\leq\deff\sigma$. On the other hand,
\[
\sup_{k,l\in\Z}|\partial g_\sigma(s^k,s^l)|=\sup_{k,l\in\Z}|\partial\sigma(k,l)|=\deff\sigma.
\]

%% file: section3.tex
\section{Homogenization and the Gromov norm}
Let $\Gamma$ be a group. There is a canonical semi-norm on bounded cohomology, the quotient semi-norm, that for a class $\alpha\in\h^n_b(\Gamma,\R)$ is given by
\[
\|\alpha\|=\inf_{f\in\alpha}\|f\|_\infty
\]
In dimension two this semi-norm is a proper norm (\cite{Iva}), called the \textit{Gromov norm}.\\\\
For a cocycle $f\in Z_\bb^n(\Gamma)$ one obviously has $\|[f]_b\|\leq\|f\|_\infty$. Using an estimate by C. Bavard we show that this is an equality in case $f$ is one of the cocycles of the previous section. For this purpose we consider the notion of a \textit{homogenous} quasi-morphim: $\varphi\in\QM(\Gamma)$ is called homogenous if $\varphi(x^n)=n\varphi(x)$ for $x\in \Gamma$, $n\in\Z$. Let $\QM^\hh(\Gamma)\subset\QM(\Gamma)$ be the subspace of homogenous quasi-morhpisms. The process of homogenizing quasi-morphisms relies on the following result:\\\\
\textbf{Lemma 3.1.} \textit{If $\{a_n\}_{n\in\N}$ is a sequence of real numbers satisfying
\[|a_n+a_m-a_{n+m}|<C\]
for $m,n\in\N$ and a constant $C>0$ then $|ma_n-a_{mn}|<(m-1)C$ and the sequence $\{\frac{a_n}{n}\}_{n\in\N}$ converges.}\\\\
\textbf{Proof.} For $m,n\in\N$ we have
\[
|ma_n-a_{mn}|=\left|\sum_{k=1}^{m-1}a_{kn}+a_n-a_{(k+1)n}\right|<(m-1)C
\]
and analogously $|na_m-a_{mn}|<(n-1)C$. I.e. $|ma_n-na_m|<(m+n-2)C$, so
\[
\left|\tfrac{a_n}{n}-\tfrac{a_m}{m}\right|<(\tfrac{1}{m}+\tfrac{1}{n}-\tfrac{2}{mn})C,
\]
which shows that $\{\frac{a_n}{n}\}$ is a Cauchy sequence.\qed
\\\\
\textbf{Proposition 3.2.} \textit{There is a direct sum decomposition
\[
\QM(\Gamma)=\QM^\hh(\Gamma)\oplus C^1_\bb(\Gamma)
\]
where the homogenous part $f^\hh$ of $f\in\QM(\Gamma)$ is given by}
\[
f^\hh(x)=\lim_{n\rightarrow\infty}\frac{f(x^n)}{n}.
\]
\textbf{Proof.} Let $f\in\QM(\Gamma)$ and $x\in\Gamma$. The previous Lemma applies to the sequence $\left\{f(x^n)\right\}_{n\in\N}$, so we can define $f^\hh$ as in the proposition. The first part of the lemma implies $|nf(x)-f(x^n)|<(n-1)C$, where $C=\deff f$, so
\[
|f(x)-f^\hh(x)|=\lim_{n\rightarrow\infty}\tfrac{1}{n}|nf(x)-f(x^n)|\leq C,
\]
which means that $f=f^\hh+f^\bb$ for some $f^\bb\in C^1_\bb(\Gamma)$. In particular, $f^\hh$ is a quasi-morphism, and it is easily seen to be homogenous. This is the unique such decomposition, since any bounded homogenous quasi-morphism is zero.\qed
\\\\
\textbf{Corollary 3.3.} \textit{For $f\in\QM(\Gamma)$ the cocycles $\partial f$ and $\partial f^\hh$ represent the same class in $\h^2_\bb(\Gamma,\R)$ and there is an isomorphism}
\[
\frac{\QM^\hh(\Gamma)}{\Hom(\Gamma,\R)}\cong\eh^2_\bb(\Gamma,\R).
\]
\textbf{Corollary 3.4.} \textit{On abelian groups all quasi-morphisms are trivial.}\\\\
\textbf{Proof.} If $\Gamma$ is abelian and $\varphi\in\QM^\hh(\Gamma)$ then
\begin{align*}
|\varphi(x)+\varphi(y)-\varphi(xy)|&=\tfrac{1}{n}|\varphi(x^n)+\varphi(y^n)-\varphi((xy)^n)|\\
&=\tfrac{1}{n}|\varphi(x^n)+\varphi(y^n)-\varphi(x^ny^n)|\rightarrow0\quad(n\rightarrow\infty).
\end{align*}
That is, $\QM^\hh(\Gamma)=\Hom(\Gamma,\R)$.\qed\\\\
\textbf{Remark.} This is a special case of a more general result: It has been shown that $\h^n_\bb(\Gamma,\R)=0$, $n\geq1$, for any \textit{amenable} group $\Gamma$, see \cite{HirThu}. By Proposition 1.1 this implies that any such group admits only trivial quasi-morphisms. Abelian groups, as well as finite groups and solvable groups are amenable.\\\\
We now determine the homogenization of the quasi-morphisms $g_\sigma\in\QM(F)$ defined in the previous section.\\\\
\textbf{Lemma 3.5.} \textit{Let $x,y\in F$.}
\begin{itemize}
\item[(i)]{\textit{If $x$ is a power then $g_\sigma^\hh(x)=0$, otherwise $g_\sigma^\hh(x)=g_\sigma(x)-\partial g_\sigma(x,x)$.}}
\item[(ii)]{\textit{If none of $x,y$ is a power then \[\partial g_\sigma^\hh(x,y)=\partial g_\sigma(x,y)+\partial g_\sigma(xy,xy)-\partial g_\sigma(x,x)-\partial g_\sigma(y,y)\]}}
\end{itemize}
\textbf{Proof.} (i) $g_\sigma^\hh(s^k)=\lim_{n\rightarrow\infty}\frac{1}{n}\sigma(kn)=0$ since the sequence $\sigma$ is bounded. If $x\in F$ is not a power then $\partial g_\sigma(x^k,x)=\partial g_\sigma(x,x)$ for $k\geq1$ and the claim follows from
\[
(n-1)\partial g_\sigma(x,x)=\sum_{k=1}^{n-1}\partial g_\sigma(x^k,x)=ng_\sigma(x)-g_\sigma(x^n).
\]
(ii) Apply the first part.\qed
\\\\
\textbf{Proposition 3.6.} $\|\partial g_\sigma^\hh\|_\infty\geq2\|\partial g_\sigma\|_\infty.$\\\\
\textbf{Proof.} Let $\varepsilon>0$, let $k,l\in\Z$ such that $d:=\partial\sigma(k,l)>\deff \sigma-\varepsilon$, and let $s,t\in S$ be distinct generators. Consider $x,y\in F$ given by
\begin{align*}
x=s^{-k}t^{-k}st^{-l}s^k,\\
y=s^lt^{-l}st^{-k}s^{-l}.
\end{align*}
These words satisfy
\[
\partial g_\sigma(x,y)=d,\quad\partial g_\sigma(x,x)=\partial g_\sigma(y,y)=\partial g_\sigma(xy,xy)=-d.
\]
Hence, by Lemma 3.5.(ii) and Proposition 1.4,
\[
\partial g_\sigma^\hh(x,y)=2d>2\,\deff\sigma-2\varepsilon=2\|\partial g_\sigma\|_\infty-2\varepsilon.
\]\qed\\
The following estimate by C. Bavard holds for any group $\Gamma$:\\\\
\textbf{Proposition 3.7.} (\cite{Bav}). \textit{Each $\varphi\in\QM^\hh(\Gamma)$ satisfies}
\[
\|[\partial\varphi]_\bb\|\geq\tfrac{1}{2}\|\partial\varphi\|_\infty.
\]
\textbf{Proposition 3.8.} $\|[\partial g_\sigma]_\bb\|=\deff\sigma$.\\\\
\textbf{Proof.} Combining the estimates of the preceeding propositions we get
\[
\|[\partial g_\sigma]_\bb\|\leq\|\partial g_\sigma\|_\infty\leq\tfrac{1}{2}\|\partial g^\hh_\sigma\|_\infty\leq\|[\partial g^\hh_\sigma]_\bb\|.
\]
Since $\partial g_\sigma$ and $\partial g^\hh_\sigma$ represent the same class these are all equalities, so $\|\partial g_\sigma]_\bb\|=\|\partial g_\sigma\|_\infty=\deff\sigma$ by Proposition 2.4.\qed\\\\
\textbf{Remark.} In particular, we have shown that the estimate in Proposition 3.7 is an equality in case of the quasi-morphisms $g_\sigma^\hh$. It has been conjectured that equality holds for any homogenous quasi-morphism on any group.\\\\ 

%% file: section4.tex
\section{Free products}
We show how the results of Section 2 can be adapted to obtain non-trivial quasi-morphisms on free products of groups.\\\\
Let $\Gamma,\Gamma'$ be groups. A map $f:\Gamma\rightarrow\Gamma'$ is called \textit{odd} if $f(x^{-1})=f(x)^{-1}$ for all $x\in\Gamma$.
We denote $\widehat{C}^1_\bb(\Gamma)\subset C^1_\bb(\Gamma)$ the space of bounded odd maps $\Gamma\rightarrow\R$.\\\\
Let $\{\Gamma_s\}_{s\in S}$, $|S|\geq2$, be a family of non-trivial groups and let $\Gamma=\ast_{s\in S}\Gamma_s$ be the associated free product. We consider the space
\[
V(\Gamma):=\prod_{s\in S}\widehat{C}^1_\bb(\Gamma_s)
\]
and its subspace of uniformly bounded families:
\[
V_0(\Gamma):=\left\{(\sigma_s)_{s\in S}\in V(\Gamma)\,\,;\,\,\sup_{s\in S}\|\sigma_s\|_\infty<\infty\right\}\subset V(\Gamma).
\]
We identify each $\Gamma_s$ with its image under the natural map $\Gamma_s\hookrightarrow\Gamma$. The \textit{factorization} of an element $x\in \Gamma$ is the unique way of writing $x$ as a product $x=x_0\cdots x_n$ such that $x_i\in \Gamma_{s_i}$ is nontrivial and $s_i\neq s_{i+1}$ for $0\leq i<n$. For $\sigma=(\sigma_s)_{s\in S}\in V_0(\Gamma)$ and $x\in\Gamma$ with factorization as above, define $g_\sigma:\Gamma\rightarrow\R$ by
\[
g_\sigma(x)=\sum_{i=0}^n\sigma_{s_i}(x_i).
\]
\textbf{Proposition 4.1.} \textit{The map $g_\sigma$ is a quasi-morphism.}\\\\
\textbf{Proof.} The argument from Proposition 2.1 holds as well in this context, so for $x,y\in\Gamma$ we have
\[|\partial g_\sigma(x,y)|\leq3\sup_{s\in S}\|\sigma_s\|<\infty.\]\\\\
\textbf{Proposition. 4.2.} \textit{The map $V_0(\Gamma)\rightarrow\eh^2_\bb(\Gamma,\R)$ given by $\sigma\mapsto[\partial g_\sigma]_b$ is a linear injection.}\\\\
\textbf{Proof.} Assume that $[\partial g_\sigma]_\bb=0$, i.e. $g_\sigma-b=\varphi$ for some $b\in C^1_\bb(\Gamma)$ and some $\varphi\in\Hom(\Gamma,\R)$. For $x\in \Gamma_s$, evaluating this equation at $x^k$ yields $\sigma_s(x^k)-b(x^k)=k\,\varphi(x)$. The left-hand side is bounded as a function of $k$, so $\varphi(x)=0$. Since $\Gamma$ is generated by the subset $\bigcup_{s\in S}\Gamma_s$, we get $\varphi=0$ and $g_\sigma$ is bounded.\\
Let $s,t\in S$ be distinct indices and let $x\in\Gamma_s, y\in\Gamma_t$. For $k\in\Z$ the equation $g_\sigma((xy^{\pm1})^k)=k(\sigma_s(x)\pm\sigma_t(y))$ holds. Since $g_\sigma$ is bounded, this implies $\sigma_s(x)\pm\sigma_t(y)=0$, so $\sigma_s(x)=\sigma_t(y)=0$ and therefore, $\sigma=0$.\qed\\\\
\textbf{Remarks}. (i) The free group $F(S)$ over a set $S$ is naturally isomorphic to the free product of a set of copies of $\Z$ indexed by $S$. Letting $\sigma_s=\sigma\in\ell^\infty$ for all $s\in S$ gives exactly the quasi-morphisms of Section 2.\\\\
(ii) The Gromov norm of the classes $[\partial g_\sigma]_\bb$ can be calculated by modifying the arguments of the previous section.\\\\
As an example we state the following\\\\
\textbf{Corollary 4.3.} \textit{For $\Gamma:=\PSL_2(\Z)$ there exists a non-trivial quasi-morphism $\Gamma\rightarrow\R$ and hence, $\dim\h^2_\bb(\Gamma,\R)\geq1$.}\\\\
\textbf{Proof.} Since $\PSL_2(\Z)\cong\Z_2\ast \Z_3$ there is a linear map
\[
V_0(\Gamma)\cong\widehat{C}^1_\bb(\Z_2)\times\widehat{C}^1_\bb(\Z_3)\hookrightarrow\eh^2_\bb(\Gamma,\R).
\]
Any odd map on $\Z_n=\{0,1,\dots,n-1\}$ is characterized by its values on the elements $1,2,\dots,\lfloor\frac{n-1}{2}\rfloor$, so we have
\[
\dim \widehat{C}^1_\bb(\Z_n)=\left\lfloor\frac{n-1}{2}\right\rfloor,
\]
and hence,
\[
\dim V_0(\Gamma)=0+1=1.
\]
\qed 

%% file: section5.tex
\section{Groups without small subgroups}
We show how the target of the quasi-morphisms of Section 2 can be replaced by a a group without small subgroups.\\\\
Let $G$ be a group with neutral element $e$, and let $d$ be a metric on $G$. For a set $X$ the distance of two maps $f,g:X\rightarrow G$ is given by \[d(f,g):=\;\sup_{x\in X}\;d(f(x),g(x)).\] If $d(f,g)<\infty$ we say $f$ and $g$ are at bounded distance. If $f$ is at bounded distance from the trivial map $x\mapsto e$ then we say $f$ is bounded and we write $\|f\|_\infty$ for this distance.\\
Let $\Gamma$ be another group. A map $f:\Gamma\rightarrow G$ is called a \textit{quasi-morphism} or an \textit{$\varepsilon$-representation} if the maps $\Gamma^2\rightarrow G$, \[(x,x')\mapsto f(xx')\quad\mbox{and}\quad(x,x')\mapsto f(x)f(x')\] are at bounded distance.\\\\
A subgroup $H\leq G$ is called $\varepsilon$-small if $H\subset B_\varepsilon(e)$. G is said to be a \textit{group without small subgroups} if there exists $\varepsilon>0$ such that every $\varepsilon$-small subgroup is trivial. For example, $\R$ equipped with the usual metric belongs to this class of groups, and the two definitions of a quasi-morphism into $\R$ agree.\\
A metric on a group is called \textit{bi-invariant} if it turns left and right translation on into isometries. Any compact Lie group carries a bi-invariant metric such that there are no small subgroups.\\\\
From now on, let $G$ be equipped with a bi-invariant metric $d$. We construct quasi-morphisms $F\rightarrow G$, where $F=F(S)$ as in Section 2.\\\\
Let $\sigma\in\ell^\infty(G)$. We define $g_\sigma:F\rightarrow G$ by $g_\sigma(x)=\sigma(k)$ if $x=s^k$ for some $s\in S$, where $\sigma$ is again extended to an odd map $\Z\rightarrow G$. In general, for $x\in F$ with factorization $x=x_0\cdots x_n$, let
\[
g_\sigma(x)=\prod_{i=0}^n g_\sigma(x_i).
\]
\textbf{Proposition 5.1.} \textit{The map $g_\sigma$ is a quasi-morphism.}\\\\
\textbf{Proof.} Let $x,y\in F$ with factorizations $x=x_0\cdots x_n$ and $y=y_0\cdots y_m$. As in Section 2 the factorization of $xy$ takes one of two possible forms. Consider the case $xy=x_0\cdots x_{n-r}\cdot z\cdot y_r\cdots y_m$, i.e. $x_{n-i}=y_i^{-1}$ for $0\leq i\leq r-2$. Due to the bi-invariance of $d$, powers cancel as in the proof of Proposition 2.1, so
\begin{align*}
&d\left(g_\sigma(xy),g_\sigma(x)g_\sigma(y)\right)\\&=d\left(\prod_{i=1}^{n-r}g_\sigma(x_i)\cdot g_\sigma(z)\cdot\prod_{i=r}^m g_\sigma(y_i)\,,\,\prod_{i=1}^n g_\sigma(x_i)\prod_{i=1}^m g_\sigma(y_i)\right)\\
&=d(g_\sigma(z),g_\sigma(x_{n-r+1})g_\sigma(y_{r-1}))\\
&\leq d(g_\sigma(z),e)+d(g_\sigma(x_{n-r+1}),e)+d(g_\sigma(y^{-1}_{r-1}),e)\leq 3\|\sigma\|_\infty.
\end{align*}
In the second case we have complete cancellation and the distance is zero.\qed\\\\
Let $G$ be a group with a bi-invariant metric $d$ and without $\varepsilon$-small subgroups, and let $\sigma$ and $g_\sigma$ be as above.\\\\
\textbf{Proposition 5.2.} \textit{Let $\varepsilon_0:=\|\sigma\|_\infty$. If $0<\varepsilon_0<\frac{\varepsilon}{2}$ then the quasi-morphism $g_\sigma$ is non-trivial, in the sense that there is no $\varphi\in\Hom(F,G)$ such that $d(g_\sigma,\varphi)\leq\varepsilon_0$.}\\\\
\textbf{Proof.} Assume there is such a $\varphi$. For $s\in S$ and $k\in\Z$ we have
\begin{align*}
d(\varphi(s^k),e)&\leq d(\varphi(s^k),g_\sigma(s^k))+d(g_\sigma(s^k),e)\leq\varepsilon_0+\varepsilon_0<\varepsilon.
\end{align*}
That is, the cyclic group $\langle\varphi(s)\rangle\leq G$ is small and therefore $\varphi(s)=e$. It follows that $\varphi$ is trivial, so $g_\sigma$ is bounded by $\frac{\varepsilon}{2}$.\\
Now pick $t\in S$ distinct from $s$ and $k,l\in\Z$. Since $(g_\sigma(s^l,t^{\pm1}))^k=g_\sigma((s^lt^{\pm1})^k)$, the cyclic group $\langle g_\sigma(s^lt^{\pm1})\rangle\leq G$ is small. So $g_\sigma(s^lt^{\pm1})=\sigma_s(l)\sigma_t(1)^{\pm1}=e$. This implies $\sigma_s(l)^{-1}=\sigma_s(l)$. Therefore the group $\langle\sigma_s(l)\rangle=\{\sigma_s(l),\sigma_s(l)^{-1}\}\leq G$ is small, so $\sigma_s(l)=e$. Since this holds for any $s\in S$ and any $l\in\Z$, we have $\|\sigma\|_\infty=0$, a contradiction.\qed\\\\
\textbf{Remark.} In the case $G=\R$ we get once more the non-triviality of the quasi-morphisms of Section 2, since $\varepsilon$-small subgroups of $\R$ are trivial for any $\varepsilon>0$. 